\documentclass{article}
\usepackage{graphicx}
\usepackage{float}
\usepackage{amsfonts}
\usepackage{amsmath}
\usepackage{amsthm}
\usepackage{amssymb}
\usepackage{mathtools}
\usepackage{commath}
\usepackage{bm}
\usepackage{color}
\usepackage{enumitem}
\usepackage{comment}
\usepackage[T1]{fontenc}
\usepackage[utf8]{inputenc}
\usepackage{multirow}
\usepackage{subcaption}
\usepackage{epsfig}

\usepackage{lineno}

\DeclareMathOperator*{\argmin}{\mathrm{argmin}}

	{\left\lbrace \begin{array}{@{} l @{} }}%
	{ \end{array}\right.}

\newtheorem{thm}{Theorem}

\newtheorem{remark}{Remark}
\newcommand{\bs}{\boldsymbol}

\newcommand{\mA}{\mathsf{A}}

\newlist{level}{itemize}{4}
\setlist[level]{label={},noitemsep,topsep=0pt}
\newcounter{algo}
\renewcommand{\thealgo}{\arabic{algo}}

\makeatletter
\def\algbackskip{\hskip-\ALG@thistlm}
\makeatother

\title{Efficient Reduced Basis Algorithm (ERBA) for kernel-based approximation}
\author{
  Francesco Marchetti$^*$, 
   Emma Perracchione$^{**}$ \\
  $^*$Dipartimento di Matematica \lq\lq Tullio Levi-Civita\rq\rq, \\
  Universit\`a di Padova, Italia \\
  \texttt{francesco.marchetti@math.unipd.it} \\
  $^{**}$CNR-SPIN, \\
  Universit\`a di Genova, Italia \\
  \texttt{perracchione@dima.unige.it} 
}
\date{May 2021}

\begin{document}

\maketitle

\begin{abstract}
    The main purpose of this work is the one of providing an efficient scheme for constructing reduced interpolation models for kernel bases. In literature such problem is mainly addressed via the well-established \emph{knot insertion} or \emph{knot removal} schemes. Such iterative strategies are usually quite demanding from a computational point of view and our goal is to study an efficient implementation for data removal approaches, namely Efficient Reduced Basis Algorithm (ERBA). Focusing on kernel-based interpolation, the algorithm makes use of two iterative rules for removing data. The former, called ERBA-$r$, is based on classical residual evaluations. The latter, namely ERBA-$p$, is independent of the function values and relies on error bounds defined by the power function. {In both cases, inspired by the so-called extended Rippa's algorithm, our ERBA takes advantage of a fast implementation.}
\end{abstract}

\section{Introduction}

We introduce the kernel-based scattered data interpolation problem following \cite{Fasshauer,Wendland05}. 
Let $ \Omega \subseteq \mathbb{R}^{d}$ and $ {\cal X} = \{  \boldsymbol{x}_i, \; i = 1,  \ldots , n\} \subset \Omega$ be a set of distinct nodes, $n\in\mathbb{N}$, where $\boldsymbol{x}_i=(x_{i,1},\dots,x_{i,d})^{\intercal}$, $i=1,\dots,n$. The scattered data interpolation problem consists in recovering an unknown function  $f: \Omega \longrightarrow \mathbb{R}$ given its values at ${\cal X}$, i.e. $\bs{f}=\bs{f}_{|{\cal X}}=(f(\bs{x}_1),\dots,f(\bs{x}_n))^{\intercal}=(f_1,\dots,f_n)^{\intercal}$. This can be achieved by imposing the interpolation conditions at ${\cal X}$. In particular, for kernel-based interpolation, the approximating function assumes the form:
\begin{equation*}\label{eq:int_problem}
    S_{f,\mathcal{X}}(\bs{x})=\sum_{i=1}^n c_i \kappa_{\varepsilon}(\bs{x},\bs{x}_i),\quad \bs{x}\in\Omega,
\end{equation*}
where $\bs{c}=(c_1,\dots,c_n)^{\intercal}\in\mathbb{R}^n$ and $\kappa_{\varepsilon}:\Omega\times\Omega\longrightarrow \mathbb{R}$ is a strictly positive definite kernel depending on a \textit{shape parameter} $\varepsilon>0$. In addition, $\kappa_{\varepsilon}$ is supposed to be \textit{radial}, i.e. there exists a univariate function $\varphi_{\varepsilon}:\mathbb{R}_{\ge0}\longrightarrow \mathbb{R}$ such that $\kappa_{\varepsilon}(\bs{x},\bs{y})=\varphi_{\varepsilon}(r)$, with $r\coloneqq\lVert \bs{x}-\bs{y}\lVert_2$, where $\lVert \cdot \lVert_2$ is the Euclidean norm. As a consequence, such a setting is commonly referred to as \emph{Radial Basis Function} (RBF) interpolation. Clarified this and with abuse of notation, we might write simply $\kappa$ instead of $\kappa_{\varepsilon}$. 
The interpolation conditions are satisfied by finding the unique vector $\bs{c}$ so that
\begin{equation}\label{eq:system}
    \mA\bs{c}=\bs{f},
\end{equation}
where $\mA=(\mA_{i,j})=\kappa(\bs{x}_i,\bs{x}_j)$, $i,j=1,\dots,n$, is the so-called interpolation (or collocation or simply kernel) matrix. The uniqueness of the solution of \eqref{eq:system} is guaranteed as long as $\kappa$ is strictly positive definite. For a more general formulation of the interpolant that involves \textit{conditionally} positive definite kernels, and for a complete overview concerning kernel-based approximation, we refer the reader e.g. to \cite{Wendland05}.

The interpolant $S_{f,\mathcal{X}}$ belongs to the space
\begin{equation*}
\mathcal{H}_{\kappa} =  \mathrm{span} \left\{ \kappa(\cdot,\boldsymbol{x}), \; \boldsymbol{x} \in \Omega \right\},
\end{equation*}
which equipped with the bilinear form $(\cdot,\cdot)_{\mathcal{H}_{\kappa}}$ is a pre-Hilbert space with reproducing kernel $\kappa$. Moreover, the completion of $\mathcal{H}_{\kappa}$ with respect to the norm $\lVert \cdot \lVert_{\mathcal{H}_{\kappa}}=\sqrt{(\cdot,\cdot)_{\mathcal{H}_{\kappa}}}$ is the so-called \textit{native space} $\mathcal{N}_{\kappa}$ associated to $\kappa$. Many upper bounds for the interpolation error, e.g. in terms of the fill-distance \cite[Theorem 14.5, p. 121]{Fasshauer} and using sampling inequalities \cite{fuselier,rieger}, are available. Here we focus on the following pointwise error bound \cite[Theorem 14.2, p. 117]{Fasshauer}: 
\begin{equation*}
    |f(\bs{x})-S_{f,\mathcal{X}}(\bs{x})|\le P_{\kappa,\mathcal{X}}(\bs{x})\lVert f \lVert_{\mathcal{N}_{\kappa}},\quad f\in \mathcal{N}_{\kappa}, \quad \bs{x}\in\Omega,
\end{equation*}
where $P_{\kappa,\mathcal{X}}$ is the so-called \textit{power function}. For our scopes, we directly define the power function as \cite[p. 116, \S 14]{Fasshauer15}:
\begin{equation}\label{eq:the_power}
    P_{\kappa,\mathcal{X}}(\bs{x})=\sqrt{\kappa(\boldsymbol{x},\boldsymbol{x})-\boldsymbol{\kappa}^{\intercal}(\boldsymbol{x}) \mA^{-1} \boldsymbol{\kappa}(\boldsymbol{x})}.
\end{equation}
where $$\bs{\kappa}(\bs{x}):=(\kappa(\bs{x},\bs{x}_1),\dots,\kappa(\bs{x},\bs{x}_n))^{\intercal}.$$ Note that \eqref{eq:the_power} splits the error into two terms. The former, i.e. the power function, only depends on the nodes and on the kernel, while the latter takes into account the function values. 

As a consequence, the power function gives information about how the interpolation error relates to the node distributions. Indeed, as for polynomial and spline interpolation, the approximation quality strongly depends on the distribution of the scattered data. In view of this, possibly starting from an initial set of nodes, many \textit{adaptive} strategies have been studied in order to construct \textit{well-behaved} interpolation designs, i.e. interpolation sets which provide an accurate reconstruction and, preferably, affordable computational costs. In particular, the so-called \textit{greedy} approaches match such purposes by iteratively adding new points to the interpolation set. The iterative rule is based on minimizing a pointwise upper bound for the interpolant. Precisely, those strategies rely on the residual of the interpolant ($f$-greedy) or on the value of the power function ($p$-greedy); refer to \cite{Dutta,Haasdonk&S,Wenzel21a,Wirtz,Wirtz1} for a general overview. These methods fall under the context of \textit{knot insertion} algorithms, which have been studied also in the context of adaptive least-squares approximation \cite[\S 21]{Fasshauer}.

Alternatively, as in our work, one could start by taking a \textit{large} interpolation set, which provides an accurate approximation, and iteratively remove nodes until the resulting interpolation error does not exceed a fixed tolerance (see e.g. \cite{Tom} for a general overview). This kind of \textit{knot removal} approaches aim to provide reduced models by neglecting as many points as possible from the initial set while preserving a target accuracy. However, they are usually computationally expensive, since many interpolants built upon large sets of nodes need to be constructed \cite{Fasshauer95}.

In this paper, in order to overcome the limitations related to the high computational complexity for reduced basis models, we take advantage of the extension of the Rippa's algorithm \cite{Rippa} recently provided in \cite{Marchetti}. More precisely, we propose the Efficient Reduced Basis Algorithm (ERBA) where the Extended Rippa's Algorithm (ERA) plays a fundamental role. Besides the residual-based rule, that leads to the ERBA-$r$, and that can be derived from the ERA, we provide an efficient scheme for computing the power function as well. The resulting ERBA-$p$ is a reduced model that takes advantage of being accurate, fast and easy to implement. 

The paper is organised as follows. In Section \ref{sec:removal_intro}, we fix some notations and we introduce the Reduced Basis Algorithm (RBA) based on knot removal strategies. Its efficient implementation, i.e. the ERBA scheme,  and theoretical results are then provided in Section \ref{sec:sezione_ciccia}, where a detailed analysis of the computational complexity of the proposed schemes is also included. Numerical tests that confirm the theoretical findings are presented in Section \ref{sec:numerics}. The conclusions are offered in Section \ref{sec:conclusions}.

%-------------------------------------------------------------------%

\section{The Reduced Basis Algorithm (RBA)}\label{sec:removal_intro}

In what follows, we present our schemes for reduced basis models based on both the residual and on the power function minimization. The latter algorithm is quasi optimal in the sense that, referring to \eqref{eq:the_power}, we only take into account the power function, neglecting the term involving the native space norm of the sought function.  This should not be seen as a drawback. Indeed, we are able to compute a quasi-optimal subset of points independently of the function values. As a consequence, it might be relevant if one has to deal with many measurements sampled at the same points. We further point out that the residual-based scheme shows strong similarities with the knot removal algorithm presented in \cite[\S 21]{Fasshauer}. We now present the RBA scheme and in the next section we focus on its fast implementation. 

Given $ {\cal X} = \{  \boldsymbol{x}_i, \; i = 1,  \ldots , n\} \subset \Omega$ be the initial set of nodes and $ {\cal X}_{s-1} = \{  \boldsymbol{x}_i, \; i = 1,  \ldots , n_{s-1}\} \subset \Omega$ be the reduced set at the $(s-1)$-th step of the algorithm ($ {\cal X}_{0} = {\cal X}$), let $\tau\in\mathbb{R}_{>0}$ be a fixed tolerance. Chosen $\rho\in\mathbb{N}$, $\rho<n$, and letting $\ell = \lfloor n_{s-1} /\rho \rfloor $, the $s$-th step of the iterative scheme, $s \geq 1$, is as follows.
\begin{enumerate}
    \item
    Partition $\mathcal{X}_{s-1}$ into $\ell$ test sets $\mathcal{X}_{s-1}^1,\dots,\mathcal{X}_{s-1}^{\ell}$, $|\mathcal{X}_{s-1}^j|=\rho^j$ with $\rho^j \in\{\rho,\dots,2\rho\}$.  Moreover, let $\overline{\mathcal{X}}_{s-1}^j\coloneqq \mathcal{X}_{s-1}\setminus\mathcal{X}_{s-1}^j$ be the training sets, $j=1,\ldots,\ell$.
    \item
    For each $\mathcal{X}_{s-1}^j$, $j=1,\dots,\ell$, compute:
    \begin{itemize}
    	\item[2a.] in the case of the residual-based scheme: 
    	\begin{equation}\label{s_int}
    	w^j = \dfrac{1}{\sqrt{\rho^j}}\lVert \bs{f}^j-\bs{S}^j \lVert_2,    
    	\end{equation}
    	where $\bs{f}^j\coloneqq\bs{f}_{|_{\mathcal{X}_{s-1}^j}}$ and $\boldsymbol{S}^j \coloneqq S_{f,\overline{\mathcal{X}}_{s-1}^j}(\mathcal{X}_{s-1}^j)$ is the evaluation vector on $\mathcal{X}_{s-1}^j$ of the interpolant $S_{f,\overline{\mathcal{X}}_{s-1}^j}$;
      	\item[2b.] in the case of the power-based scheme:
      	\begin{equation*}\label{p_int}
      	w^j = \dfrac{1}{\sqrt{\rho^j}} \lVert \bs{P}^j\lVert_2,
      \end{equation*}	
       where $\bs{P}^j\coloneqq P_{\kappa,\overline{\mathcal{X}}_{s-1}^j}(\mathcal{X}_{s-1}^j)$ is the evaluation vector on $\mathcal{X}_{s-1}^j$ of the power function $P_{\kappa,\overline{\mathcal{X}}_{s-1}^j}$.
    \end{itemize}
    \item
    Choose
    $$ j^{\star}=\argmin_{j\in\{1,\dots,\ell\}}{w^j}. $$
    \item 
    Let $r_{s-1} = w^{j^{\star}}$:
    \begin{itemize}
        \item[4a.]
        if $r_{s-1}\le\tau$, define $\mathcal{X}_{s} = \mathcal{X}_{s-1}\setminus \mathcal{X}_{s-1}^{j^{\star}}$ and proceed iteratively with the $s$-th step; 
        \item[4b.]
        if $r_{s-1}>\tau$, $\mathcal{X}_{s} = \mathcal{X}_{s-1}$ is the final interpolation set obtained by the procedure.
    \end{itemize}
\end{enumerate}

%The convergence of the residual-based scheme has been studied in many context and for many basis functions, see e.g. \cite{Tom} for a general overview. As far as the convergence of the proposed error-based scheme is concerned,  let 
%$$
%V(\mathcal{X}_{s}) \subset V(\mathcal{X}_{s-1}) \subset \ldots \subset V(\mathcal{X}_{1}) \subset {\cal N}_{\kappa}, 
%$$
%be linear subspaces of the nested sets $\mathcal{X}_i$, $i=1,\ldots,s$. 
%Then, by using the Newton basis \cite{SS}, we can consider an  orthonormal basis for $V({\cal X}_s)$, namely $\{b_k\}_{k=1}^{m_1}$, $m_1=|{\cal X}_s|.$ With these ingredients we can write \cite[Lemma 5]{Haasdonk&S}:
%$$
%P_{\kappa,V({\cal X}_s)} (\boldsymbol{x}) = \lVert \kappa(\cdot,\boldsymbol{x})-S_{\kappa(\cdot,\boldsymbol{x}),V({\cal X}_s)}(\boldsymbol{x}) \lVert_{{\cal N}_{\kappa}},
%$$
%and hence, letting $m_2=|{\cal X}_{s-1}|,$ 
%\begin{equation*}
%\begin{split}
%      P^2_{\kappa, V({\cal X}_{s-1})}  (\boldsymbol{x}) & =  \kappa(\boldsymbol{x},\boldsymbol{x}) - \sum_{k=1}^{m_2} (b_k(\boldsymbol{x}))^2, \\
%         & =   
%\kappa(\boldsymbol{x},\boldsymbol{x}) - \sum_{k=1}^{m_1} (b_k(\boldsymbol{x}))^2 - \sum_{k=m_1+1}^{m_2} (b_k(\boldsymbol{x}))^2, \\
%& = P^2_{\kappa, V({\cal X}_s)}  (\boldsymbol{x})  - \sum_{k=m_1+1}^{m_2} (b_k(\boldsymbol{x}))^2 \leq  P^2_{\kappa, V({\cal X}_s)}   (\boldsymbol{x}).
%\end{split}
%\end{equation*}
%This implies that:
%$$
%\lVert P_{\kappa,V({\cal X}_{s-1})}  \lVert_{L_\infty(\Omega)} \leq \lVert P_{\kappa,V({\cal X}_{s})}  \lVert_{L_\infty(\Omega)}.
%$$

The convergence of the residual-based scheme has been studied in many context and for many basis functions, see e.g. \cite{Tom} for a general overview. As far as the convergence of the proposed error-based scheme is concerned, 
it is ensured by the fact that \cite[Lemma 5]{Haasdonk&S}: 
$$
\lVert P_{\kappa,V({\cal X}_{s-1})}  \lVert_{L_\infty(\Omega)} \leq \lVert P_{\kappa,V({\cal X}_{s})}  \lVert_{L_\infty(\Omega)},
$$
being 
$
V(\mathcal{X}_{s}) \subset V(\mathcal{X}_{s-1}), 
$ linear subspaces of the nested sets $\mathcal{X}_{s}$ and $\mathcal{X}_{s-1}$.

From now on, in analyzing the computational complexity, we assume without loss of generality that $\rho$ divides $n_s$ at each step of the algorithm, and thus $\rho_j\equiv\rho$ holds true. Therefore, at each step, a classical implementation of this method requires to solving $n_s$ different $(n_s-\rho)\times (n_s-\rho)$ linear systems, in the residual-based case, or the inversion of a $(n_s-\rho)\times (n_s-\rho)$ matrix in the power-based setting. Hence, in both situations, each step is computationally demanding. In \cite[\S 21]{Fasshauer}, a similar residual-based approach has been speeded up by using Rippa's algorithm in the case $\rho=1$. In the next section, inspired by ERA \cite{Marchetti} we provide a fast implementation of the proposed algorithm.

%-------------------------------------------------------------------%

\section{Efficient Reduced Basis Algorithm (ERBA)}\label{sec:sezione_ciccia}

To study the computational complexity of the proposed schemes, we have to focus on the calculation of the residuals (or of the power function) that is performed at each step. Then, fixed a certain iteration $s=1,2,\dots$, let us introduce the following notations: $\mathcal{X}_s\coloneqq \mathcal{X}$, $n_s\coloneqq n$, and let $\bs{p}\coloneqq (p_1,\dots,p_{\rho})^{\intercal}$ be the vector of $\rho$ indices related to the elements of the subset $\mathcal{X}_s^j\coloneqq\mathcal{V}$, $|\mathcal{X}_s^j|=\rho$. We also adopt the following notations:
\begin{gather*}
\bs{f}_{\bs{p}}\coloneqq (f_{i})_{i\in \bs{p}},\quad \bs{f}^{\bs{p}}\coloneqq (f_{i})_{i\notin \bs{p}},\\
\mA^{\bs{p},\bs{p}}\coloneqq (\mA_{i,j})_{i,j\notin \bs{p}},\quad \mA_{\bs{p},\bs{p}}\coloneqq (\mA_{i,j})_{i,j\in \bs{p}},\quad \mA_{\bs{p},:}\coloneqq (\mA_{i,j})_{i\in \bs{p},j\in\{1,\dots,n\}},\\
\bs{\kappa}(\bs{x}):=(\kappa(\bs{x},\bs{x}_1),\dots,\kappa(\bs{x},\bs{x}_n))^{\intercal}, \\
\bs{k}_{\bs{p}}(\bs{x})= (\bs{\kappa}(\bs{x})_{i})_{i\in \bs{p}},\quad \bs{k}^{\bs{p}}(\bs{x})= (\bs{\kappa}(\bs{x})_{i})_{i\notin \bs{p}},
\end{gather*}
In the following we explain our efficient strategy for the implementation of both the residual and power-based schemes.

\subsection{The ERBA-$r$}

In step 2a. for the computation of $\boldsymbol{S}^j$ in equation \eqref{s_int}, we need to compute the interpolant and hence solve a system that leads to a matrix inversion.
Precisely, let us consider
\begin{equation*}
S_{f,{\overline{\cal V}}}(\bs{x})\coloneqq S_{f,{\cal X}}^{(\bs{p})}(\bs{x})=\sum_{i\notin \bs{p}} c^{(\bs{p})}_i \kappa_{\varepsilon}(\bs{x},\bs{x}_i),\quad \bs{x}\in\Omega,\; \bs{x}_i\in\mathcal{X},
\end{equation*}
where $\bs{c}^{(\bs{p})}\coloneqq \big(c^{(\bs{p})}_{i}\big)_{i\notin \bs{p}}$ is determined by solving
\begin{equation*}
\mA^{\bs{p},\bs{p}}\bs{c}^{(\bs{p})}=\bs{f}^{\bs{p}}.
\end{equation*}
We are interested in computing the residual
\begin{equation}\label{eq:residual}
\bs{e}_{\bs{p}} \coloneqq \bs{f}_{\bs{p}}-((\bs{c}^{(\bs{p})})^{\intercal}\bs{k}^{\bs{p}}(\mathcal{V}))^{\intercal}
\end{equation}
being $\bs{k}^{\bs{p}}(\mathcal{V})=(\kappa(\bs{x}_i,\bs{x}_j))_{i\in\bs{p},j\notin \bs{p}}$ a $(n-\rho)\times \rho$ matrix.\\
Supposing $n/\rho=\ell\in\mathbb{N}$, a classic approach would lead to the resolution of $\ell$ different $(n-\rho)\times(n-\rho)$ linear systems, and thus, since usually $n \gg \rho$, to a computational cost of about $\mathcal{O}(n^4)$. In case $\rho=1$ one could use the Rippa's algorithm reducing the complexity to $\mathcal{O}(n^3)$. In case of more folds, i.e. $\rho >1$, we take advantage of the ERA scheme \cite{Marchetti} for which  \eqref{eq:residual} can be computed as
\begin{equation}\label{eq:rippa_residual}
\bs{e}_{\bs{p}} = (\mA^{-1}_{\bs{p},\bs{p}})^{-1}\bs{c}_{\bs{p}}.
\end{equation}
Doing in this way, we need to invert an $n\times n$ matrix and then to solve $\ell$ different $\rho \times \rho$ linear systems. Therefore, the computational cost is about $\mathcal{O}(n^3)$, leading to a significant saving (cf. \cite[\S 2.2]{Marchetti}).

%-------------------------------------------------------------------%

\subsection{The ERBA-$p$}

Here, we focus on the computation of the power function vector (see \eqref{eq:the_power} and step 2b. of the algorithm presented in Section \ref{sec:removal_intro})
\begin{equation}\label{eq:power_vector}
\bs{P}^{\mathcal{V}}=\sqrt{\mathrm{diag}(\bs{k}_{\bs{p}}(\mathcal{V})-\bs{k}^{\bs{p}}(\mathcal{V})^{\intercal}(\mA^{\bs{p},\bs{p}})^{-1}\bs{k}^{\bs{p}}(\mathcal{V}))},
\end{equation}
where $\mathrm{diag}(\cdot)$ is the diagonal operator. The vector $\bs{P}^{\mathcal{V}}$ is usually computed by means of a for-loop over the elements of $\mathcal{V}$ (see \cite[Program 17.1]{Fasshauer}). Hence,  supposing again $n/\rho=\ell\in\mathbb{N}$, we need to invert $\ell$ matrices of dimension $(n-\rho)\times (n-\rho)$, leading to a computational cost of about $\mathcal{O}(\ell(n-\rho)^3)=\mathcal{O}(n(n-\rho)^3/\rho)$. Since $n \gg \rho$, the cost for each step of the algorithm is approximately $\mathcal{O}(n^4)$.

To provide a fast computation of the power function vector $\bs{P}^{\mathcal{V}}$, we first recall that, letting $\mathsf{R} \in \mathbb{R}^{n_1\times n_2}$ and $\mathsf{S} \in \mathbb{R}^{n_2\times n_1}$, we have that
\begin{equation}\label{eq:hada}
\mathrm{diag}(\mathsf{R}\mathsf{S})=\mathrm{sum}(\mathsf{R}\odot\mathsf{S}^{\intercal}), 
\end{equation}
where $\odot$ denotes the Hadamard (or pointwise) product. 
Then, we obtain the following result.
\begin{thm}\label{thm:power_to_rippa}
	The vector \eqref{eq:power_vector} can be computed as
	\begin{equation*}
	\bs{P}^{\mathcal{V}}=\sqrt{\mathrm{sum}(((\mA^{-1}_{\bs{p},\bs{p}})^{-1}\mA^{-1}_{\bs{p},:})\odot \bs{k}(\mathcal{V})^{\intercal})},
	\end{equation*}
	where $\mathrm{sum}(\cdot)$ denoted the sum-by-rows operator.
\end{thm}
\begin{proof}
	Putting together \eqref{eq:residual} and \eqref{eq:rippa_residual}, we have that
	\begin{align*}
	(\mA^{-1}_{\bs{p},\bs{p}})^{-1}\bs{c}_{\bs{p}}=\bs{f}_{\bs{p}}-\bs{k}^{\bs{p}}(\mathcal{V})^{\intercal}(\mA^{\bs{p},\bs{p}})^{-1}\bs{f}^{\bs{p}},\\
	\bs{k}^{\bs{p}}(\mathcal{V})^{\intercal}(\mA^{\bs{p},\bs{p}})^{-1}\bs{f}^{\bs{p}}=\bs{f}_{\bs{p}}-(\mA^{-1}_{\bs{p},\bs{p}})^{-1}\bs{c}_{\bs{p}}.
	\end{align*}
	Since $\bs{f}$ is arbitrary if the function $f$ is not fixed, we can substitute it with $\bs{k}^{\bs{p}}(\mathcal{V})$ in the equation. Then,
	\begin{align*}
	\bs{k}^{\bs{p}}(\mathcal{V})^{\intercal}(\mA^{\bs{p},\bs{p}})^{-1}\bs{k}^{\bs{p}}(\mathcal{V})=\bs{k}_{\bs{p}}(\mathcal{V})-(\mA^{-1}_{\bs{p},\bs{p}})^{-1}\mA^{-1}_{\bs{p},:}\bs{k}(\mathcal{V}).
	\end{align*}
	Hence, recalling \eqref{eq:power_vector}, we get
	\begin{equation*}
	\bs{P}^{\mathcal{V}}=\sqrt{\mathrm{diag}((\mA^{-1}_{\bs{p},\bs{p}})^{-1}\mA^{-1}_{\bs{p},:}\bs{k}(\mathcal{V}))},
	\end{equation*}
	Finally, by applying \eqref{eq:hada} to $(\mA^{-1}_{\bs{p},\bs{p}})^{-1}\mA^{-1}_{\bs{p},:}$ and $\bs{k}(\mathcal{V})$, we conclude the proof.
\end{proof}
By adopting the scheme proposed in Theorem \ref{thm:power_to_rippa}, the computational cost at each step is about $\mathcal{O}(n^3)$, indeed we need to invert a unique $n\times n$ matrix and to perform other \textit{minor} calculations that are negligible as long as $n \gg \rho$. The proposed strategy is therefore faster than the classical framework, as confirmed by the experiments carried out in the next section.

%-------------------------------------------------------------------%

\section{Numerics}\label{sec:numerics}

In what follows, we perform some numerical experiments to prove the efficiency of the proposed ERBA. The tests have been carried out in \textsc{Matlab} on a Intel(R) Core(TM) i7-1165G7 CPU@2.80GHz processor. The software is available for the scientific community at
\begin{equation*}
\texttt{https://github.com/cesc14/ERBA}\:.
\end{equation*}
Let $\Omega=[-1,1]^2$ and let $f:\Omega\longrightarrow\mathbb{R}$ be a function defined as
\begin{equation*}
    f(\bs{x})=\frac{1}{1+(x_1-0.5)^2+(x_2+0.2)^2},\quad \bs{x}=(x_1,x_2).
\end{equation*}
For the tests, we take the strictly positive definite kernel
\begin{equation*}
\varphi(r)= e^{-\varepsilon r} ,\quad\textrm{Matérn $C^0$},
\end{equation*}
and we set $\varepsilon = 1$. 

\begin{figure}
\begin{center}
\begin{subfigure}[b]{0.49\textwidth}
\includegraphics[width=\textwidth]{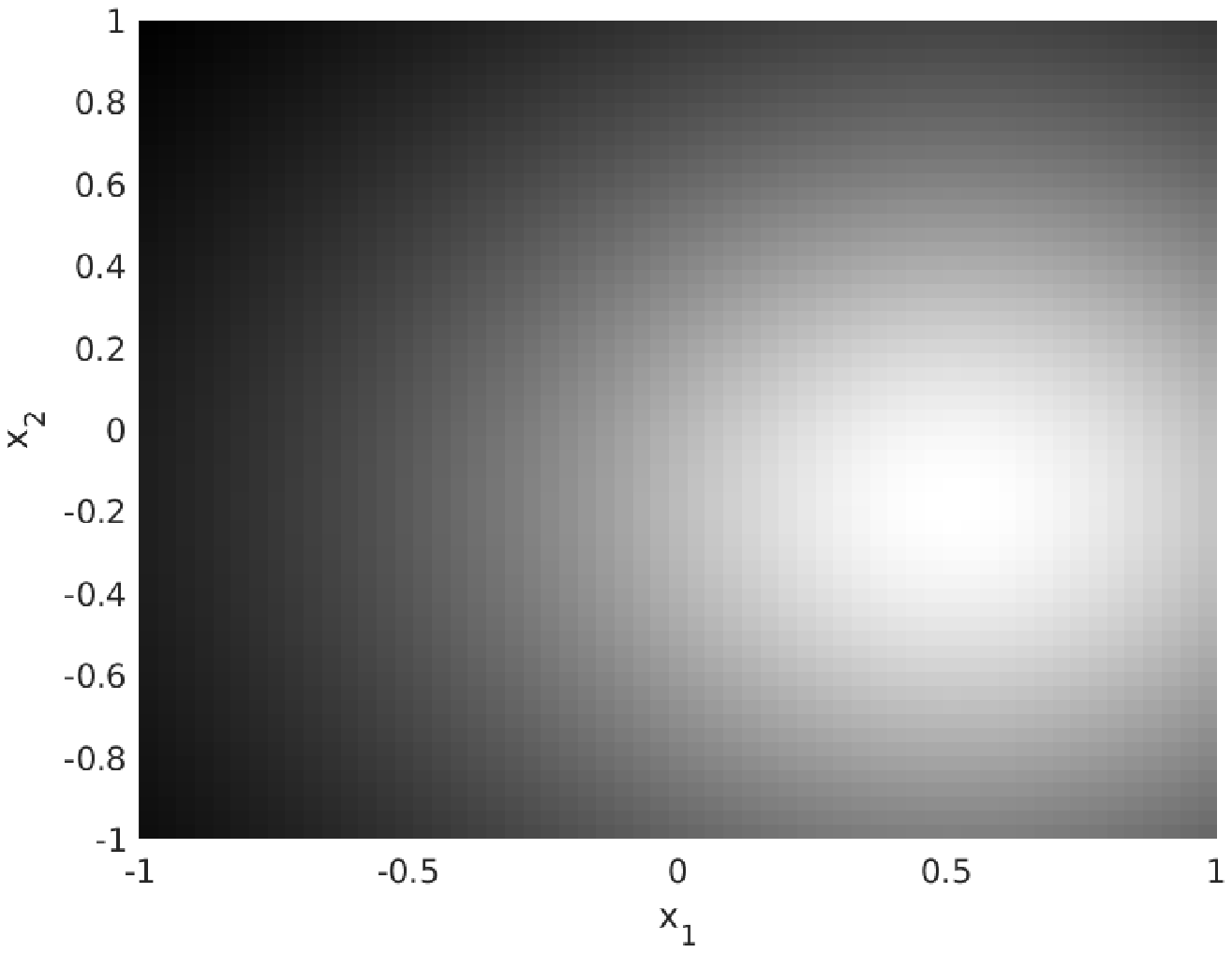}
\caption{}\label{fig:1_1}
\end{subfigure}
\begin{subfigure}[b]{0.49\textwidth}
\includegraphics[width=\textwidth]{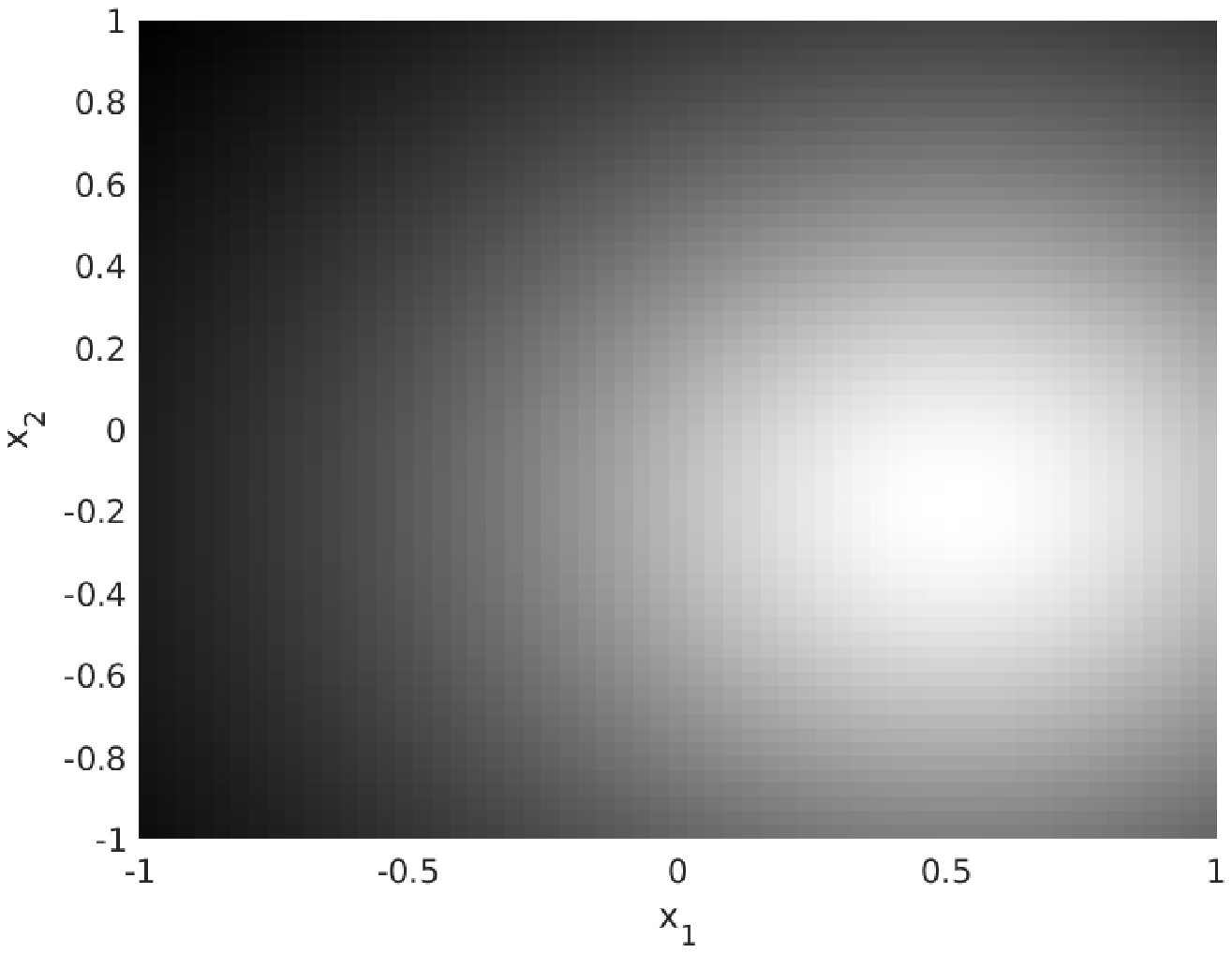}
\caption{}\label{fig:1_2}
\end{subfigure}
\end{center}
\caption{(a): the function $f$. (b) the interpolant constructed via the data set $\mathcal{X}$.}
\label{fig:1}
\end{figure}

As interpolation data set $\mathcal{X} \subset \Omega$, we consider a $n\times n$ grid, with $n=25$. Moreover, we take a $m\times m$ evaluation grid $\Xi$ with $m=60$. The associated RMSE computed on $\Xi$ is $e_{\mathcal{X}}=9.69{\rm E}-05$. 
In Figure \ref{fig:1}, we plot the function $f$ and the interpolant $S_{f,\mathcal{X}}$ evaluated on $\Xi$. 

\subsection{Testing ERBA-$r$}

In the following, we take $\rho=3$ and we set $\tau = 2  {e}_{\mathcal{X}}$.

\begin{figure}
\begin{center}
\begin{subfigure}[b]{0.47\textwidth}
\includegraphics[width=\textwidth]{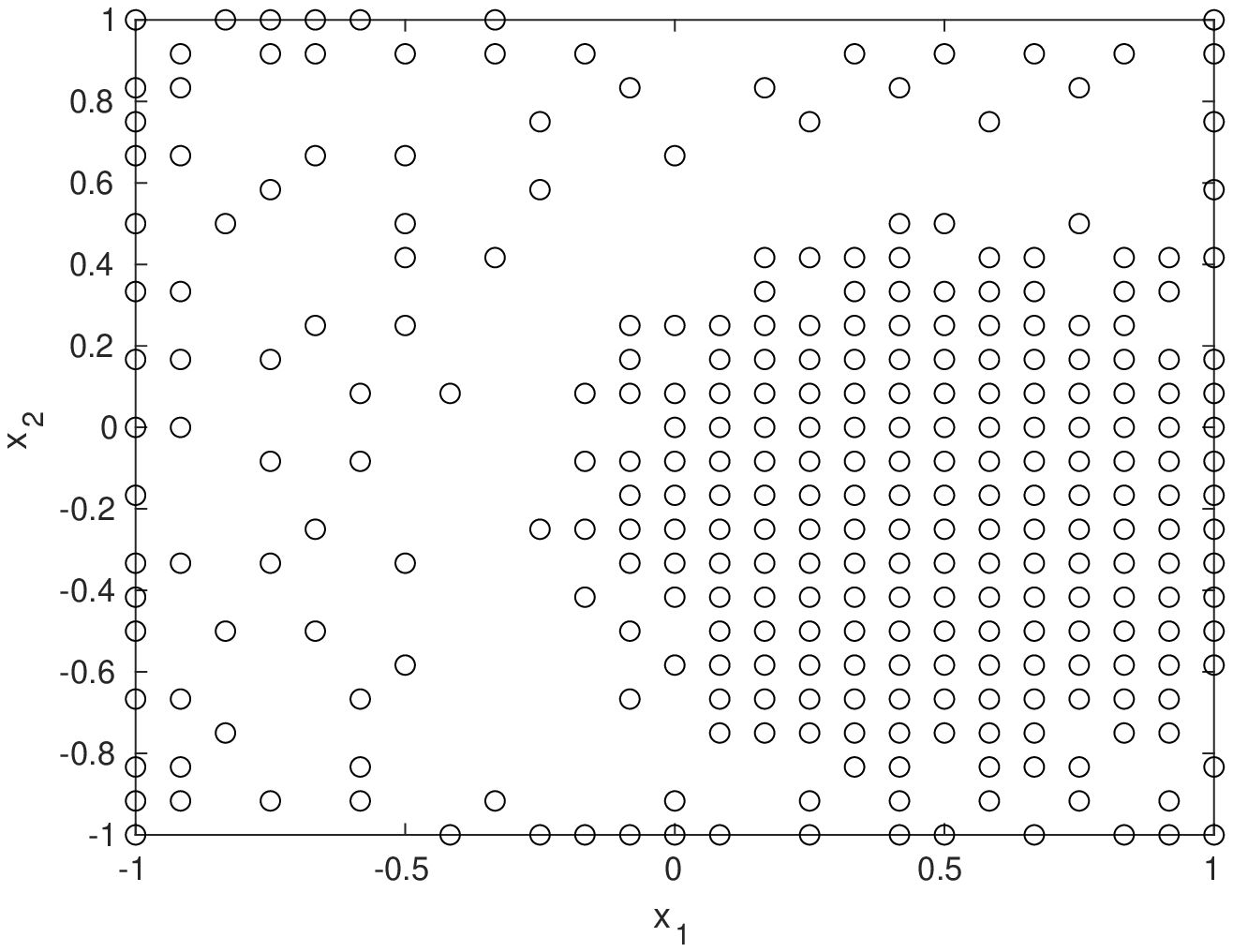}
\caption{}\label{fig:2_1}
\end{subfigure}
\begin{subfigure}[b]{0.47\textwidth}
\includegraphics[width=\textwidth]{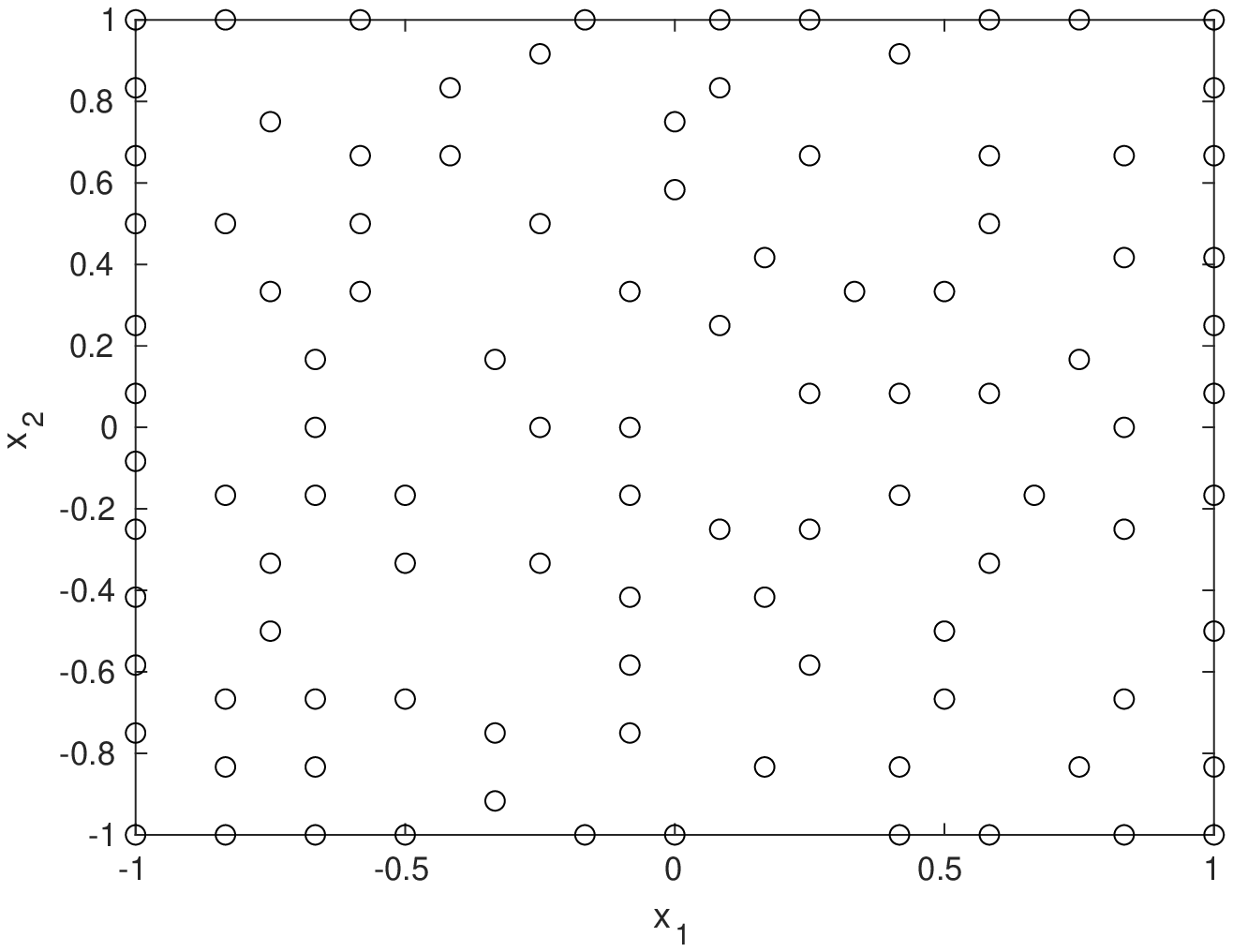}
\caption{}\label{fig:2_2}
\end{subfigure}
\end{center}
\caption{(a): the reduced ERBA-$r$ data set. (b): the reduced ERBA-$p$ data set.  }
\label{fig:2}
\end{figure}

We apply ERBA-$r$ and we denote as $\mathcal{X}_s \subset\mathcal{X}$ the set of resulting reduced interpolation data which are depicted in Figure \ref{fig:2_1}, while in Figure \ref{fig:3_1} we show the interpolant constructed with such nodes. The efficiency of the computational strategies proposed in Section \ref{sec:sezione_ciccia} is compared to the classical implementation in Table \ref{table:2}, where some further details concerning the experiments are also reported.

\begin{table}
\centering
\begin{tabular}{|| c | c | c | c ||} 
 \hline
 $|\mathcal{X}_s|$ & ${e}_{\mathcal{X}_s}$ & ERBA CPU times & RBA CPU times \\ [0.5ex] 
 \hline
$298$ & $1.29{\rm E}-04$ & $3.14{\rm E}+00$ & $5.30{\rm E}+01$\\ [0.5ex]
\hline
\end{tabular}
\caption{Results for the ERBA-$r$ scheme. CPU times are in seconds.}
\label{table:2}
\end{table}

\subsection{Testing ERBA-$p$}

Here, we set $\rho=3$ and the tolerance $\tau = 2\lVert \bs{P}^{\Xi} \lVert_{2}/m$, being $\bs{P}^{\Xi}$ the power function vector constructed via the nodes $\mathcal{X}$ and evaluated on $\Xi$. 

\begin{figure}
\begin{center}
\begin{subfigure}[b]{0.49\textwidth}
\includegraphics[width=\textwidth]{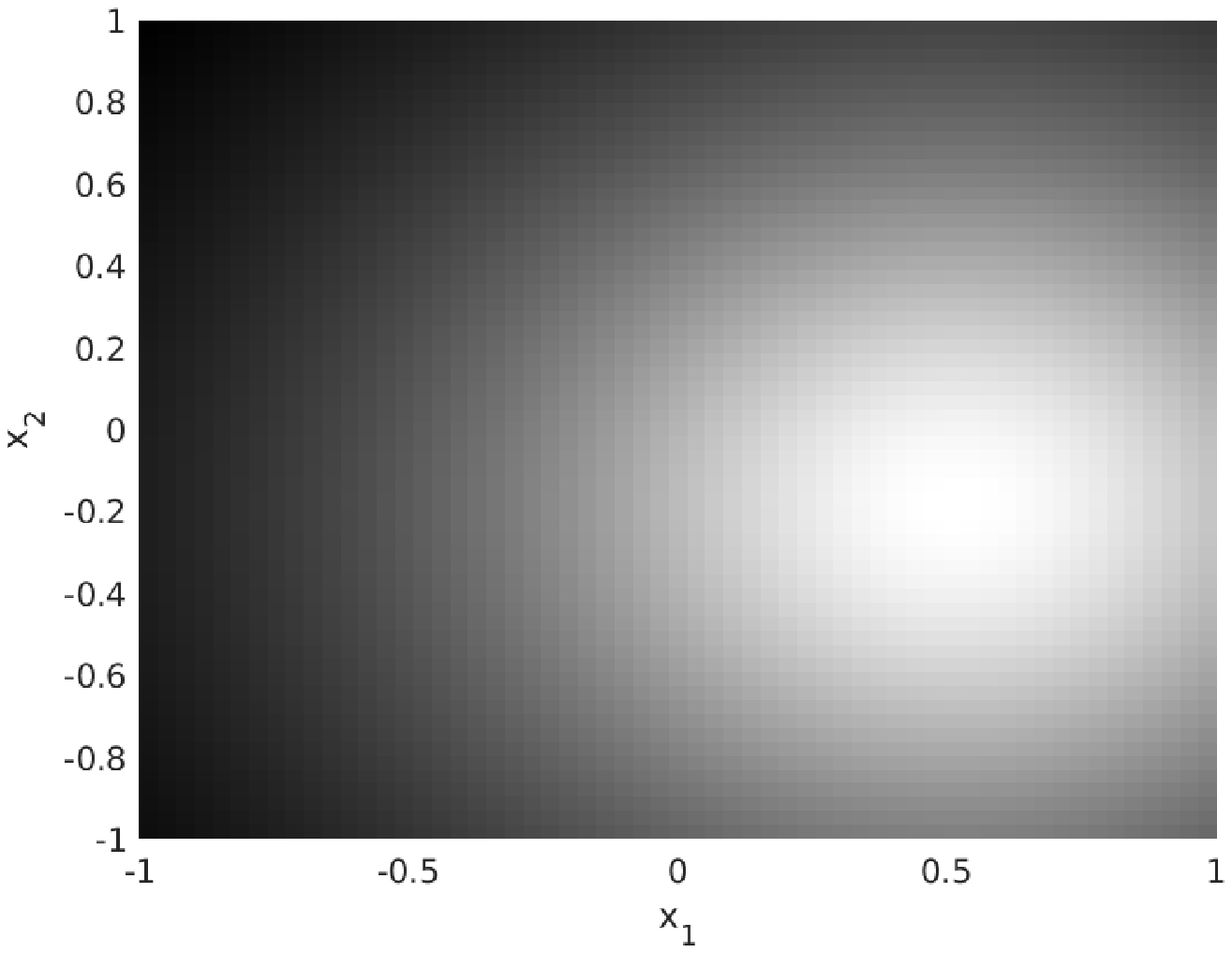}
\caption{}\label{fig:3_1}
\end{subfigure}
\begin{subfigure}[b]{0.49\textwidth}
\includegraphics[width=\textwidth]{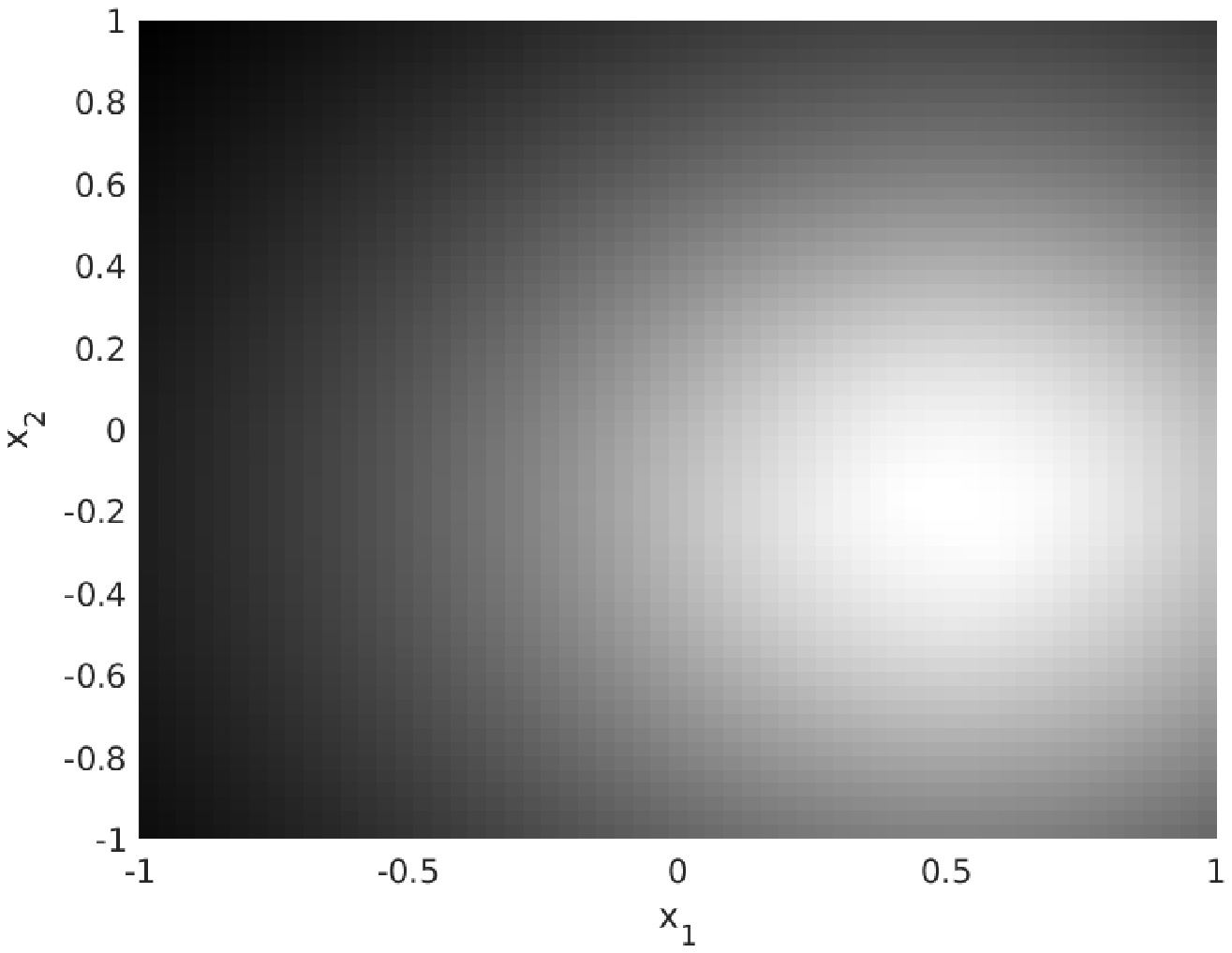}
\caption{}\label{fig:3_2}
\end{subfigure}
\end{center}
\caption{(a): the interpolant constructed via the reduced ERBA-$r$ data set. (b): the interpolant constructed via the reduced ERBA-$p$ data set.  }
\label{fig:3}
\end{figure}

As previously done, we display the results obtained via ERBA-$p$ in Figures \ref{fig:2_2} and \ref{fig:3_2} and in Table \ref{table:1}.

\begin{table}
\centering
\begin{tabular}{|| c | c | c | c ||} 
 \hline
 $|\mathcal{X}_s|$ & ${e}_{\mathcal{X}_s}$ & ERBA CPU times & RBA CPU times \\ [0.5ex] 
 \hline
$103$ & $2.41{\rm E}-03$ & $4.73{\rm E}+00$ & $5.28{\rm E}+02$\\ [0.5ex]
\hline
\end{tabular}
\caption{Results for the ERBA-$p$ scheme. CPU times are in seconds.}
\label{table:1}
\end{table}

\subsection{A focus on computing times}

\begin{figure}
\begin{center}
\begin{subfigure}[b]{0.47\textwidth}
\includegraphics[width=\textwidth]{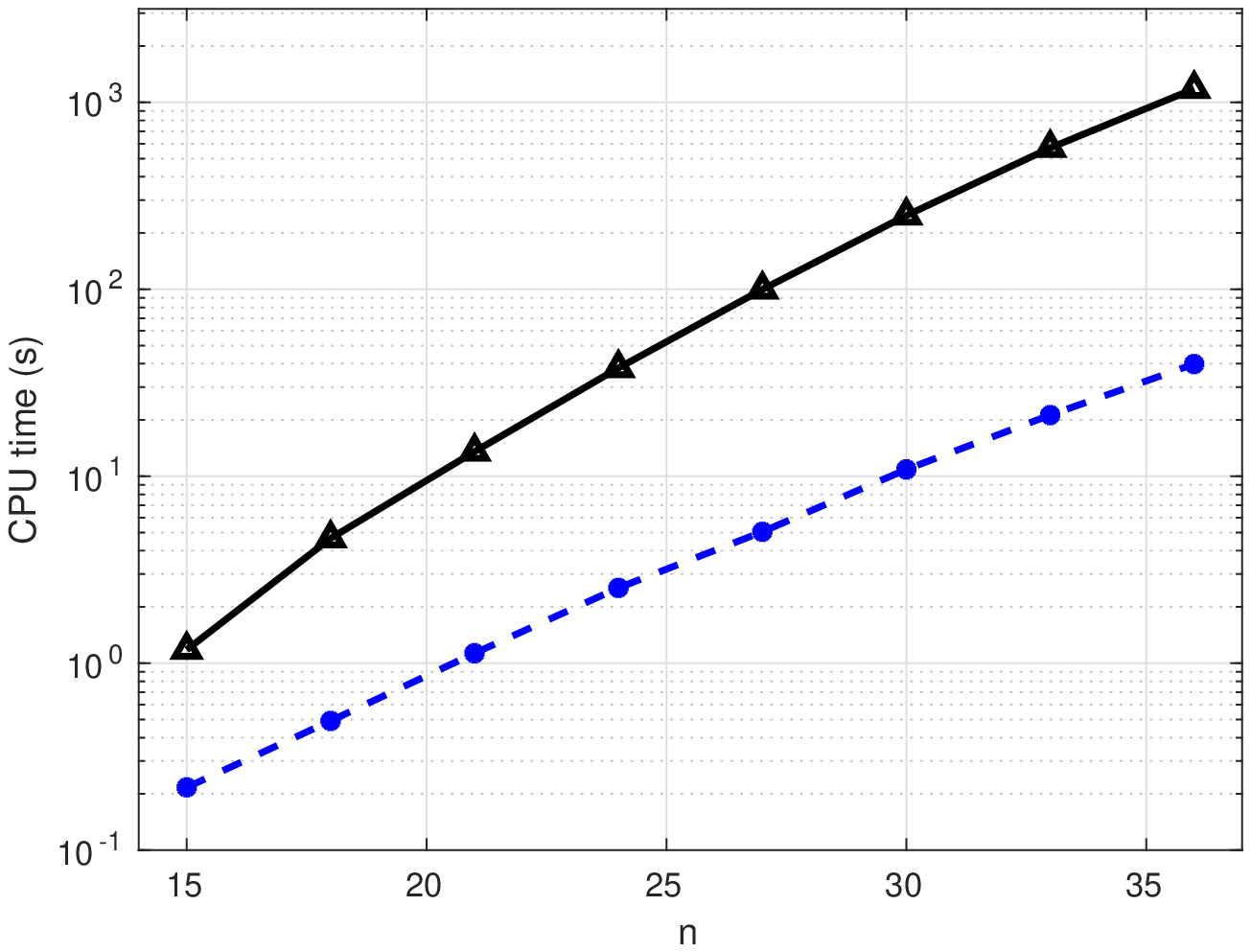}
\caption{}\label{fig:4_1}
\end{subfigure}
\begin{subfigure}[b]{0.47\textwidth}
\includegraphics[width=\textwidth]{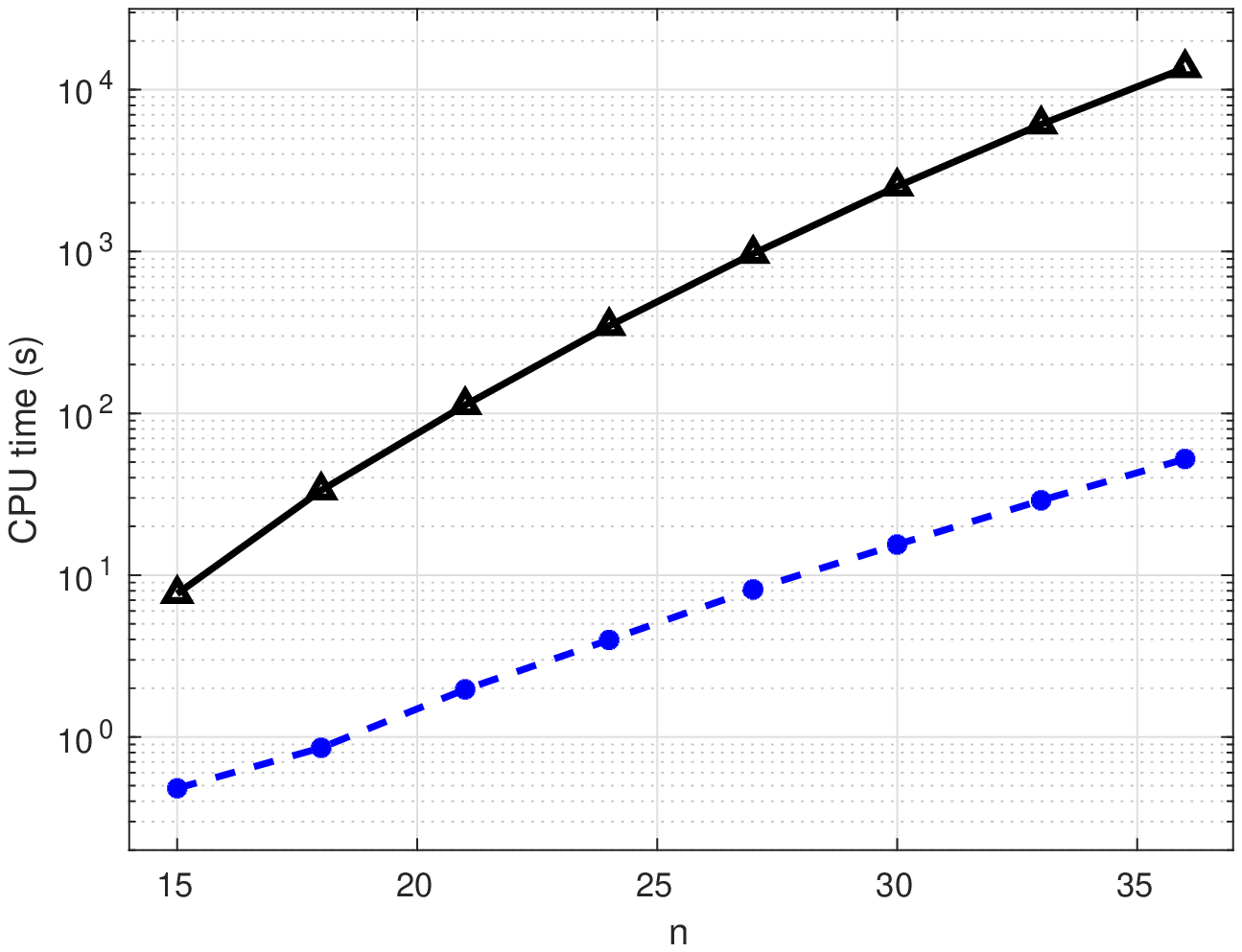}
\caption{}\label{fig:4_2}
\end{subfigure}
\end{center}
\caption{The CPU times required by ERBA (dashed blue line) and RBA (solid black line). (a): the residual based case. (b): the power based case.}
\label{fig:4}
\end{figure}

As a final experiment, we compare the CPU times of the algorithms by varying the grid size of the initial data set $\mathcal{X}$. Precisely, we take $n=15+3k$, $k=0,\ldots,7$. The results are reported in Figure \ref{fig:4}. 

\section{Discussion and conclusions}\label{sec:conclusions}

We have investigated a fast computation of a knot removal scheme. We have implemented two different strategies: the first one in based on classical residual-based schemes while the second one relies on power function error bounds. The latter is independent on the function values and tends to return \emph{quasi-uniform} data (cf. \cite{SDW} and \cite[Theorem 15 \& 19 \& 20]{Wenzel21}), while the former, as expected, keeps points where the corresponding function has steep gradients. In both cases we are able to significantly speed up the algorithm thanks to our implementation. 

Work in progress consists in extending the proposed tool to other bases, e.g. splines \cite{Campagna2019} and to variably scaled kernels \cite{Campi21}. Moreover, it can be helpful for validating the shape parameter in RBF interpolation; refer e.g. to \cite{Cavoretto}.

\section*{Acknowledgements}
This research has been accomplished within Rete ITaliana di Approssimazione (RITA) and the UMI Group TAA Approximation Theory and Applications; partially funded by GNCS-IN$\delta$AM and  by  the ASI - INAF grant  \lq\lq Artificial Intelligence for the analysis of solar FLARES data (AI-FLARES)\rq\rq.

\end{document}